\documentclass[12pt, reqno]{amsart}
\usepackage{amsmath, amsthm, amscd, amsfonts, amssymb, graphicx, color}
\usepackage[bookmarksnumbered, colorlinks, plainpages]{hyperref}
\hypersetup{colorlinks=true,linkcolor=red, anchorcolor=green, citecolor=cyan, urlcolor=red, filecolor=magenta, pdftoolbar=true}

\newtheorem{theorem}{Theorem}[section]
\newtheorem{lemma}[theorem]{Lemma}

\theoremstyle{definition}

\theoremstyle{remark}

\numberwithin{equation}{section}

\begin{document}
\setcounter{page}{1}

\title[Derivations, local and 2-local derivations]{Derivations, local and 2-local derivations of standard operator algebras}

\author[J. He, H. Zhao, G. An]{Jun He$^{1}$, Haixia Zhao$^{2}$ and Guangyu An$^{*}$}

\address{$^{1,2}$ Department of Mathematics, Anhui Polytechnic University,
Wuhu 241000, China.}
\email{\textcolor[rgb]{0.00,0.00,0.84}{hejun$_{-}$12@163.com}; \textcolor[rgb]{0.00,0.00,0.84}{1755808875@qq.com}}

\address{$^{3*}$ Department of Mathematics, Shaanxi University of Science and Technology,
Xi'an 710021, China.}
\email{\textcolor[rgb]{0.00,0.00,0.84}{anguangyu310@163.com}}



\date{Received: xxxxxx; Revised: yyyyyy; Accepted: zzzzzz.
\newline \indent $^{*}$ Corresponding author}

\begin{abstract}
Let $X$ be a Banach space over field $\mathbb F$ ($\mathbb R$ or $\mathbb C)$.
Denote by $B(X)$ the set of all bounded linear operators
on $X$ and by $F(X)$ the set of all finite rank operators on $X$.
A subalgebra $\mathcal A$ of $B(X)$ is called a standard operator algebra if
$\mathcal A$ contain $F(X)$. We give a brief proof of a well-known result that every derivation
from $\mathcal A$ into $B(X)$ is inner. There is another classical result that every local derivation on $B(X)$ is a derivation.
We extend the result by proving that every local derivation from $\mathcal A$ into $B(X)$ is a derivation. Based on these two results,
we prove that every 2-local derivation from $\mathcal A$ into $B(X)$ is a derivation.
\end{abstract}\maketitle

\section{Introduction}
Let $\mathcal{A}$ be an algebra and $\mathcal M$ be an $\mathcal A$-bimodule.
Suppose that $\delta$ is a linear mapping from $\mathcal A$ into $\mathcal M$.
$\delta$ is called a \emph{derivation}
if $\delta(xy)=\delta(x)y+x\delta(y)$ for each $x, y$ in $\mathcal A$,
and $\delta$ is called an \emph{inner derivation}
if there exist an element $m$ in $\mathcal M$ such that $\delta(x)=mx-xm$.
Clearly, every inner derivation is a derivation. But the converse is not so trivial.
It is a classical problem to identify those algebras on which every derivation is an inner derivation.
In \cite{Kadison2, sakai}, R. Kadison and S. Sakai independently prove that every derivation
on a von Neumann algebra is an inner derivation.
In \cite{Chernoff}, P. Chernoff proves that every derivation from a standard operator algebra
into $B(X)$ is an inner derivation. In this paper, we shall give a brief proof of this result.

In 1990, R. Kadison \cite{Kadison}, D. Larson and A. Sourour \cite{Larson} independently introduce
the concept of local derivation.
A linear mapping $\delta$ from $\mathcal{A}$ into $\mathcal{M}$ is called a \emph{local derivation}
if for every $x$ in $\mathcal{A}$, there exists a derivation $\delta_x$ (depends on $x$) from $\mathcal{A}$ into $\mathcal{M}$,
such that $\delta(x)=\delta_x(x)$.
In \cite{Kadison}, R. Kadison proves that every continuous local
derivation from a von Neumann algebra into its dual Banach
module is a derivation. In \cite{Larson}, D. Larson and A. Sourour prove that every local derivation
from $B(X)$ into itself is a derivation.
For more information about local derivations, we refer to \cite{Johnson,Crist,Li1,Li2,Li3,Zhu,zhujun}.
In this paper, we shall extend the result of D. Larson and A. Sourour \cite{Larson} to that
every local derivation from a standard operator algebra into $B(X)$ is a derivation.

The concept of 2-local derivation is firstly introduced by P. $\check{S}$emrl \cite{Semrl} in 1997.
A mapping (not necessarily linear) $\delta$ from $\mathcal{A}$ into $\mathcal{M}$ is called a \emph{2-local derivation}
if for each $x,y$ in $\mathcal{A}$, there exist a derivation
$\delta_{x,y}$ (depends on $x,y$) from $\mathcal{A}$ into $\mathcal{M}$, such that  $\delta(x)=\delta_{x,y}(x)$ and
$\delta(y)=\delta_{x,y}(y)$. P. $\check{S}$emrl \cite{Semrl} proves that every 2-local derivation
from $B(H)$ into itself is a derivation for a separable Hilbert space $H$.
Many other authors study 2 local derivations and there are several important results, we refer to \cite{Ayupov1,Ayupov2,Ayupov3,Zhang,Kim,Me}.
In this paper, we prove that
every 2-local derivation from a standard operator algebra into $B(X)$ is a derivation.

\section{Main results}\

We shall firstly review some simple properties of $B(X)$, especially about finite rank operators.
Throughout this section, $X$ is a Banach space over field $\mathbb F$ ($\mathbb R$ or $\mathbb C)$,
and $X^*$ is the set of all bounded linear functionals on $X$.
Denote by $B(X)$ the set of all bounded linear operators on $X$.
An operator $A\in B(X)$ is said to be of finite rank if the range of $A$ is a finite dimensional subspace of $X$.
Denote by $F(X)$ the set of all finite rank operators on $X$.
$\mathcal A$ is a standard operator algebra, which means that $\mathcal A$ is a subalgebra of $B(X)$ and $\mathcal A\supseteq F(X)$.
For each $x$ in $X$ and $f$ in $X^{*}$, one can define an operator $x\otimes f$ by $(x\otimes f)y=f(y)x$ for all $y$ in $X$.
Obviously, $x\otimes f\in B(X)$. If both $x$ and $f$ are nonzero, then $x\otimes f$ is of rank one.

The following properties are evident and will be used frequently in the proof behind.
For each $x,y$ in $X$ , $f,g$ in $X^{*}$, and $A,B$ in $B(X)$,\\
(1) $(x\otimes f)A=x\otimes(fA)$, and $A(x\otimes f)=(Ax)\otimes f$;\\
(2) $(x\otimes f)(y\otimes g)=f(y)(x\otimes g)$;\\
(3) every operator in $F(X)$ of rank $n$ can be written into the form $\sum_{i=1}^{n}x_i\otimes f_i$, where $x_i\in X$ and $f_i\in X^*$;\\
(4) $F(X)$ is a two-side ideal of $B(X)$;\\
(5) $F(X)$ is a separating set of $B(X)$, which means that $AF(X)=0$ implies $A=0$ and $F(X)A=0$ implies $A=0$.

Now we are in position to give our main results. The following Theorem \ref{2} can be found in \cite{Chernoff}.
But we shall give another brief proof here.
\begin{lemma}\label{1}
Every derivation $\delta$ from standard operator algebra $\mathcal A$ into $B(X)$ is continuous.
\end{lemma}
\begin{proof}
Assume that $\{T_n\}\subseteq\mathcal A$ is a sequence converging to $0$,
and $\{\delta(T_n)\}$ converges to $T$.
According to the closed graph Theorem, to prove $\delta$ is continuous,
it is sufficient to show that $T=0$.

For each $x,y$ in $X$ and $f,g$ in $X^*$, we have that,
$$\delta(x\otimes fT_ny\otimes g)=\delta(f(T_ny)x\otimes g)=f(T_ny)\delta(x\otimes g).$$
By $T_n\rightarrow0$, it follows that $\delta(x\otimes fT_ny\otimes g)\rightarrow0$.
On the other hand, we have that
$$\delta(x\otimes fT_ny\otimes g)=\delta(x\otimes f)T_ny\otimes g+x\otimes f\delta(T_n)y\otimes g+x\otimes fT_n\delta(y\otimes g).$$
By $T_n\rightarrow0$  and $\delta(T_n)\rightarrow T$,
it follows that
$$\delta(x\otimes fT_ny\otimes g)\rightarrow x\otimes fTy\otimes g=f(Ty)x\otimes g.$$
It implies that $f(Ty)x\otimes g=0$ for each $x,y\in X$ and $f,g\in X^*$. Hence we can obtain that $T=0$.
The proof is complete.
\end{proof}

\begin{theorem}\label{2}
Every derivation $\delta$ from standard operator algebra $\mathcal A$ into $B(X)$ is an inner derivation.
\end{theorem}
\begin{proof}
Let $x_0$ be in $X$ and $f_0$ be in $X^*$ such that $f_0(x_0)=1$.
Define a mapping $T$ from $X$ into itself by
$$Tx=\delta(x\otimes f_0)x_0$$
for all $x$ in $X$.
By Lemma \ref{1}, we know that $\delta$ is continuous.
So it is not difficult to check that $T$ is also continuous.
Obviously $T$ is linear. Hence we have $T\in B(X)$.

For each $x$ in $X$ and $A$ in $\mathcal A$,
we have that
$$TAx=\delta(Ax\otimes f_0)x_0=\delta(A)x\otimes f_0x_0+A\delta(x\otimes f_0)x_0=\delta(A)x+ATx.$$
It implies that $\delta(A)=TA-AT$ for all $A$ in $\mathcal A$.
That is to say, $\delta$ is an inner derivation. The proof is complete.
\end{proof}

\begin{theorem}\label{3}
Every local derivation $\delta$ from standard operator algebra $\mathcal A$ into $B(X)$ is a derivation.
\end{theorem}

\begin{proof}
We shall firstly consider the case that $\mathcal A$ contains unit $I$.

By the definition of local derivation,
there exists a derivation $\delta_I$ such that $\delta(I)=\delta_I(I)=0$.
Let $A$, $B$ and $C$ be in $\mathcal A$ with $AB=BC=0$.
There exists a derivation $\delta_B$ such that $\delta_B(B)=\delta(B)$.
Thus we have that
\begin{align}
A\delta(B)C=A\delta_B(B)C=\delta_B(ABC)-\delta_B(A)BC-AB\delta_B(C)=0.                               \label{201}
\end{align}

Let $A_0$ and $B_0$ be in $\mathcal A$ with $A_0B_0=0$.
Define a bilinear mapping $\phi_1$ from $\mathcal A\times\mathcal A$ into $B(X)$ by
$$\phi_1(X,Y)=X\delta(YA_0)B_0$$
for each $X$ and $Y$ in $\mathcal A$. By \eqref{201}, we see that
$XY=0$ implies $\phi_1(X,Y)=0$. Let $A$ be in $\mathcal A$ and $R$ be in $F(X)$.
By \cite[Theorem 4.1]{M. Bresar 3}, we have that
$$\phi_1(R,A)=\phi_1(RA,I),$$
that is,
$$R\delta(AA_0)B_0=RA\delta(A_0)B_0.$$
Since $F(X)$ is a separating set of $B(X)$,
we obtain that $\delta(AA_0)B_0=A\delta(A_0)B_0$ whenever $A_0B_0=0$.

Define a bilinear mapping $\phi_2$ from $\mathcal A\times\mathcal A$ into $B(X)$ by
$$\phi_2(X,Y)=\delta(AX)Y-A\delta(X)Y$$
for each $X$ and $Y$ in $\mathcal A$.
By the previous discussion, we see that $XY=0$ implies $\phi_2(X,Y)=0$.
Let $B$ be in $\mathcal A$ and $R$ be in $F(X)$.
Again by \cite[Theorem 4.1]{M. Bresar 3},
we have that
$$\phi_2(B,R)=\phi_2(I,BR).$$
By $\delta(I)=0$, it implies that
$$\delta(AB)R-A\delta(B)R=\delta(A)BR.$$
Since $F(X)$ is a separating set of $B(X)$,
we obtain that
$$\delta(AB)=A\delta(B)+\delta(A)B$$
for each $A,B$ in $\mathcal A$.
That is to say, $\delta$ is a derivation.

Next we shall consider the case that $\mathcal A$ do not contain unit $I$.
In this case, denote the unital algebra $\mathcal A\oplus\mathbb{F}I$ by $\widetilde{\mathcal A}$.

For every linear mapping $\phi$ from $\mathcal A$ into $B(X)$,
we can define a linear mapping $\widetilde{\phi}$ from $\widetilde{\mathcal A}$ into $B(X)$
by $\widetilde{\phi}(A+\lambda I)=\phi(A)$ for all $A\in\mathcal A $ and $\lambda\in\mathbb F$.
For all $A,B\in\mathcal A $ and $\lambda,\mu\in\mathbb F$, we have
$$\widetilde{\phi}((A+\lambda I)(B+\mu I))=\widetilde{\phi}(AB+\lambda B+\mu A+\lambda\mu I)=\phi(AB)+\lambda\phi(B)+\mu\phi(A)$$
and
$$\widetilde{\phi}(A+\lambda I)(B+\mu I)+(A+\lambda I)\widetilde{\phi}(B+\mu I)=\phi(A)B+A\phi(B)+\lambda\phi(B)+\mu\phi(A).$$
It implies that $\phi$ is a derivation if and only if $\widetilde{\phi}$ is a derivation.

Since $\delta$ is a local derivation from $\mathcal A$ into $B(X)$.
For each $A\in\mathcal A $ and $\lambda\in\mathbb F$, there exists a derivation $\delta_A$ such that $\delta(A)=\delta_A(A)$.
Moreover, we have $\widetilde{\delta}(A+\lambda I)=\delta(A)=\delta_A(A)=\widetilde{\delta_A}(A+\lambda I)$.
It means that $\widetilde{\delta}$ is a local derivation from $\widetilde{\mathcal A}$ into $B(X)$.
By the result of the case that $\mathcal A$ contains unit, $\widetilde{\delta}$ is a derivation.
Hence $\delta$ is also a derivation. The proof is complete.
\end{proof}

\begin{theorem}
Every 2-local derivation $\delta$ from standard operator algebra $\mathcal A$ into $B(X)$ is a derivation.
\end{theorem}

\begin{proof}
For all $A\in\mathcal A$ and $S\in F(x)$,
by the definition of 2-local derivation,
there exists a derivation $d$ from $\mathcal A$ into $B(X)$ such that
$\delta(A)=d(A)$ and $\delta(S)=d(S)$.
By Theorem \ref{2}, $d$ is an inner derivation.
Thus we have
$$\delta(A)S+A\delta(S)=d(A)S+Ad(S)=d(AS)=AST-TAS$$
for some $T\in B(X)$.
It follows that $tr(\delta(A)S+A\delta(S))=0$, where $tr$ is the trace mapping on $F(X)$.
Moreover, $tr(\delta(A)S)=-tr(A\delta(S))$.

Now, for each $A,B\in\mathcal A$ and $S\in F(X)$,
we have
\begin{align*}
tr(\delta(A+B)S)&=-tr((A+B)\delta(S))=-tr(A\delta(S))-tr(B\delta(S))\\
&=tr(\delta(A)S)+tr(\delta(B)S)\\
&=tr((\delta(A)+\delta(B))S).
\end{align*}

Let $C=\delta(A+B)-\delta(A)-\delta(B)$,
we obtain $tr(CS)=0$.
By taking $S=x\otimes f$, where $x\in X$ and $f\in X^*$ are chosen arbitrily,
we have $tr(Cx\otimes f)=f(Cx)=0.$
It follows that $C=0$, i.e. $\delta(A+B)=\delta(A)+\delta(B)$.
That is to say $\delta$ is additive.
In addition, by the definition of 2-local derivation, it is easy to see that $\delta$ is homogeneous.
Hence $\delta$ is linear, moreover, a local derivation.
Then by Theorem \ref{3}, $\delta$ is a derivation.
The proof is complete.
\end{proof}

\bibliographystyle{amsplain}

\begin{thebibliography}{99}
\bibitem{Ayupov1} S. Ayupov, K. Kudaybergenov, 2-Loacl derivations and automorphisms on $B(H)$, J. Math. Anal. Appl., 395(2012), 15-18.

\bibitem{Ayupov2} S. Ayupov, K. Kudaybergenov, 2-local derivations on von Neumann algebras, Positivity, 19 (2014), 445-455.

\bibitem{Ayupov3} S. Ayupov, K. Kudaybergenov,  A. Peralta, A survey on local and 2-local derivations on C*- and von Neuman algebras, Topics in functional analysis and algebra, 73-126, Contemp Math., 672, 2016.

\bibitem{M. Bresar 3} M. Bre\v{s}ar, Multiplication algebra and maps determined by zero products.
Linear Multilinear Algebra, 60 (2012), 763--768.

\bibitem{Chernoff} P. Chernoff, Representations, automorphism and derivation of some operator algebras, J. Funct. Anal., 12 (1973), 275-289.

\bibitem{Christensen} E. Christensen, Derivations of nest algebras, Math. Ann., 229(1977), 155-161.

\bibitem{Crist} R. Crist, Local derivations on operator algebras, J. Funct. Anal., 135 (1996), 72-92.

\bibitem{Me} J. He, J. Li, G. An, W. Huang. Characterizations of 2-local derivations and local Lie derivations on some algebras.
Sib. Math. J., 59(2018), 721-730.

\bibitem{Kim} S. Kim, J. Kim, Local automorphisms and derivations on $M_n(\mathbb{C})$, Proc. Amer. Math. Soc., 132(2004), 1389-1392.

\bibitem{Li1} D. Hadwin, J. Li, Local derivations and local automorphisms, J. Math. Anal. Appl., 290(2004), 702-714.

\bibitem{Li2} D. Hadwin, J. Li, Local derivations and local automorphisms on some algebras, J. Operator Theory, 60(2008), 29-44.

\bibitem{Johnson} B. Johnson, Local derivations on $C^*$-algebras are derivations, Trans. Amer. Math. Soc., 353(2001), 313-325.

\bibitem{Kadison2} R. Kadison, Derivations of operator algebras, Ann. Math., 83(1966), 280-293.

\bibitem{Kadison} R. Kadison, Local derivations, J. Algebra, 130(1990), 494-509.

\bibitem{Larson} D. Larson, A. Sourour, Local derivations and local automorphisms, Proc. Sympos. Pure Math., 51(1990), 187-194.

\bibitem{Li3} J. Li, Z. Pan, Annihilator-preserving maps, multipliers and local derivations, Linear Algebra Appl., 432(2010), 5-13.

\bibitem{sakai} S. Sakai, Derivations of $W^{*}$-algebras, Ann. Math., 83(1966), 273-279.

\bibitem{Semrl} P. $\check{S}$emrl, Local automorphisms and derivations on $B(H)$, Proc. Amer. Math. Soc., 125(1997), 2677-2680.

\bibitem{Zhang} J. Zhang, H. Li, 2-Loacl derivations on digraph algebras, Acta Math. Sinica(Chin. Ser.), 49(2006), 1401-1406.

\bibitem{Zhu} J. Zhu, Local derivation of nest algebras, Proc. Amer. Math. Soc., (123)1995, 739-742.

\bibitem{zhujun} J. Zhu, C. Xiong, Bilocal derivations of standard operator algebras, Proc. Amer. Math. Soc., 125(1997), 1367-1370.


\end{thebibliography}

\end{document}